\documentclass[12pt]{article}

\usepackage{verbatim}

\makeatletter

   \@addtoreset{equation}{section}
\makeatother

\begin{document}

\vskip 12mm

\begin{center} 
{\Large \bf  An Analytic Formula for Numbers of Restricted Partitions from Conformal Field Theory}
\vskip 10mm
{ \large  Dimitri Polyakov$^{a,b,c}$\footnote{email:polyakov@scu.edu.cn;polyakov@sogang.ac.kr;
twistorstring@gmail.com}

\vskip 8mm
$^{a}$ {\it  Center for Theoretical Physics, College of Physical Science and Technology}\\
{\it  Sichuan University, Chengdu 6100064, China}\\
\vskip 2mm
$^{b}${\it  Max-Planck-Instutut fuhr Gravitationsphysik (Albert-Einstein-Institut), Am Muhlenberg 1, D-14476 Potsdam, Germany}

\vskip 2mm

$^{c}$ {\it Institute of Information Transmission Problems (IITP)}\\
{\it  Bolshoi Karetny per. 19/1, Moscow 127994, Russia}\\

}
\end{center}

\vskip 15mm

\begin{abstract}

We study the correlators of irregular vertex operators in two-dimensional
conformal field theory (CFT) in order to propose an exact analytic formula
for calculating  numbers of partitions, that is:

1) for given $N,k$, finding the total number $\lambda(N|k)$
of length $k$ partitions of $N$: $N=n_1+...+n_k;0<n_1\leq{n_2}...\leq{n_k}$.

2) finding the total number $\lambda(N)=\sum_{k=1}^N\lambda(N|k)$ of partitions
of a natural number $N$

We propose an exact analytic expression for $\lambda(N|k)$
by relating two-point short-distance correlation functions of irregular
vertex operators in $c=1$ conformal field theory ( the form 
of the operators is established in this paper):
with the first correlator counting the
partitions in the upper half-plane and the second one obtained
from  the first correlator by conformal transformations of the form
$f(z)=h(z)e^{-{i\over{z}}}$ where $h(z)$ is regular and non-vanishing
at $z=0$.

The final formula for $\lambda(N|k)$ is given in terms of regularized
($\epsilon$-ordered)
finite series in the higher-derivative Schwarzians and incomplete Bell polynomials
of the above conformal transformation at $z=i\epsilon$ ($\epsilon\rightarrow{0}$)

\end{abstract}

\vskip 12mm

\setcounter{footnote}{0}

\section{\bf Introduction} 
Let
\begin{eqnarray}
N=n_1+n_2+...+n_k (1{\leq{k}}\leq{N})
\nonumber \\
0<n_1{\leq}n_2...{\leq}n_k
\end{eqnarray}
be the length $k$ partition of a natural number $N$,
$\lambda(N|k)$ be the number of such length $k$ partitions of $N$ 
and 
\begin{eqnarray}
\lambda(N)=\sum_{k=1}^N\lambda(N|k)
\end{eqnarray}
be the total number of partitions.
Physically, $\lambda(N|k)$  and $\lambda(N)$
count the number of Young diagrams with $N$ cells 
and $k$ rows , and  the total number $\lambda(N)$
of the diagrams with $N$ cells, and therefore are related
to counting irreducible representations for 
higher-spin fields with spin value $N$.
As it is well-known from number theory,
obtaining exact analytic expressions
for $\lambda(N)$  and especially for $\lambda(N|k)$
(say, in terms of some finite series) is a
hard long-standing problem.
For $\lambda(N)$, various asymptotic formulae are known
for the large $N$ limit. The oldest and perhaps the best-known
formula for $\lambda(N)$ was obtained by Ramanujan and Hardy
in 1918 \cite{HR}
and is given by:
\begin{eqnarray}
\lambda(N)\sim{1\over{4N{\sqrt{3}}}}e^{\pi{\sqrt{{2N\over3}}}}
\nonumber
\end{eqnarray}
There are several improvements of this formula,
notably by Rademacher \cite{RH}, \cite{erdos} who expressed $\lambda(N)$
in terms of infinite convergent series:

\begin{eqnarray}
\lambda(N)={1\over{\pi{\sqrt{2}}}}\sum_{n=1}^\infty{\sqrt{n}}\alpha_n(N)
{{d}\over{dN}}\lbrace
{{sinh\lbrack{{\pi}\over{n}}
{\sqrt{{2\over3}(N-{1\over{24}})}}\rbrack}\over{{\sqrt{N-{1\over{24}}}}}}\rbrace
\nonumber
\end{eqnarray}

where
{
\begin{eqnarray}
\alpha_n(N)=\sum_{{0\leq{m}\leq{n}};{\lbrack}m|n{\rbrack}}
e^{i\pi(s(m,n)-{{2Nm}\over{n}})}
\nonumber
\end{eqnarray}}
with the notation ${\lbrack}m|n{\rbrack}$ implying the sum over $m$ taken 
over the values of $m\
$ relatively prime to $n$
and
{
\begin{eqnarray}
s(m,n)={1\over{4n}}\sum_{k=1}^{n-1}
cot({{\pi{k}}\over{n}})cot({{\pi{km}}
\over{n}})
\nonumber
\end{eqnarray}}
is the Dedekind sum for co-prime numbers.

The problem of finding $\lambda(N|k)$ is well-known to be even more
tedious (see e.g. \cite{almkvistf, almkvists, suresh}
for the discussion of Ramanujan-Rademacher type asymptotics for 
the restricted partitions).
In this paper, we study the two-point short-distance correlator
of irregular vertex operators
in Conformal Field Theory \cite{BPZ, DF}
that counts the number $\lambda(N|k)$ of restricted partitions , reproducing
the well-known generating function for the partitions, when computed
in the upper half-plane. 
One of these operators is the special case
of rank one irregular vertex operators \cite{AGT, GT_2012, G_2009, RP},
 that
can be physically interpreted as a ``dipole''
in the Liouville theory (in the same sense that regular vertex operators,
or primary fields,
are the ``charges'');
 another is related to
a class of analytic solutions in open string field theory \cite{wit, selff},
interpolating between flat and AdS backgrounds.

Next, we investigate the behaviour of this
correlator under the peculiar class of  conformal transformations
that shrink the dipole's size to zero and reduce the correlator to 
contribution from zero modes of the irregular vertices.
This leads to nontrivial identities involving the 
 restricted partitions,  expressing them in terms of generalized higher-derivative 
Schwarzians of these conformal transformations.
In particular, this allows to express the number $\lambda(N|k)$ of the
restricted partitions in terms of the finite series of the generalized 
Schwarzians and incomplete Bell polynomials of the conformal transformations considered, 
leading to the main result of this work.
Taking the short-distance limit in the correlation 
functions is necessary in order to be able to
integrate the Ward identities, accounting for the 
non-global part of the two-dimensional conformal symmetry (or physically, the
``spontaneous  breaking''  of
the conformal symmetry for transformations with
non-zero Schwarzians, considered in our work, i.e. other than fractional-linear).
In general, such an integration is hard to perform and the correlators, computed
in different coordinates (related by the transformation) differ by the infinite
sum over  Schwarzians and their higher-derivative counterparts. 
This difference, however, becomes controllable
in the short-distance limit and, for the conformal transformations with the asymptotics,
considered in this paper, can be compensated by a relatively simple factor,
derived in our work.

\section{ \bf Generating Function for Partitions: the Correlator} 

As it is well-known, $\lambda(N)$ and $\lambda(N|k)$ can be realized as 
expansion coefficients of the following (respectively) generating functions:

\begin{eqnarray}
F(x)=\prod_{n=1}^\infty{1\over{1-x^n}}=\sum_{N=0}^\infty\lambda(N)x^N
\nonumber \\
F(x,y)=
\prod_{n=1}^\infty{1\over{1-yx^n}}
=\sum_{N=0,k=0;k{\leq{N}}}^\infty\lambda(N|k)x^Ny^k
\nonumber
\end{eqnarray}

Unfortunately, these generating functions by themselves are 
not very helpful for elucidating explicit expressions for the partition numbers: taking their
derivatives just gives trivial identities
of the form $\lambda(N|k)=\lambda(N|k)$.
For this reason, our strategy in this work  will be to

1) identify the two-point correlators in Conformal Field Theory ($CFT$) counting the partitions
(reproducing the generating function $F(x,y)$)  in certain coordinates, namely,
an upper half-plane

2)using the conformal symmetry and  suitable conformal transformations
(identified below), derive the identities for the generating function, casting
it in terms of an expression, making it possible to
obtain an exact analytic expression $\lambda(N|k)$ in terms 
of the finite series.
For simplicity, in this paper we shall concentrate on $c=1$ CFT
(free massless bosons in two dimensions). 
With some effort, it is straightforward to identify the two-point correlator,
counting the partitions on the upper half-plane.
This correlator is given by

\begin{eqnarray}
G(\alpha,\beta,\epsilon)=
<U_{\alpha}(z_1)V_\beta(z_2)>|_{z_1=i\epsilon;z_2=0}
\nonumber \\
\epsilon>0
\end{eqnarray}

where

\begin{eqnarray}
U_\alpha(z)=
:\prod_{n=0}^\infty{1\over{1-{{\alpha^n\partial^n\phi}\over{n!}}}}:(z)
\nonumber \\
V_\beta(w)=:e^{\beta\partial\phi}:(w)
\nonumber
\end{eqnarray}

where the :: symbol stands for the normal ordering 
of operators in two-dimensional CFT, $\phi$ is $D=2$ boson 
(e.g. a Liouville field, an open string's
target space coordinate or a bosonized ghost),
$\epsilon\rightarrow{0}$, $\alpha$ and $\beta$  are the parameters that are 
introduced to
control the generating function for the partitions.
In this work, both $U_\alpha$ and $V_\beta$ are understood in terms
of formal series in $\alpha$ and $\beta$, with each term in the series being normally 
ordered by definition. 
To simplify notations , here and below we shall often use
the partial derivative symbol for $z$-derivatives, even though in
our case it coincides with ordinary derivative, since we only consider 
holomorphic sector. 

Indeed, expanding in $\alpha$:
{
\begin{eqnarray}
U_\alpha=\sum_{k=1}^\infty\sum_{n_1\leq...{\leq}n_k=0}^\infty
\sum_{p_1,...,p_k}
{{\alpha^{p_1n_1+p_2n_2+...+p_kn_k}}\over{(n_1!)^{p_1}...(n_k!)^{p_k}}}
\nonumber \\
\times
:(\partial^{n_1}\phi)^{p_1}
...(\partial^{n_k}\phi)^{p_k}:
\nonumber
\end{eqnarray}}
 using the operator product expansion (OPE):
\begin{eqnarray}
\partial^n\phi(z):e^{\beta\partial\phi}:(w)
\sim{{(-1)^{n+1}n!{\beta}:e^{\beta\partial\phi}:(w)}\over{(z-w)^{n+1}}}
+...
\nonumber
\end{eqnarray}
and introducing
{$N={\sum_{}}p_kn_k$} one easily calculates

\begin{eqnarray}
G(\alpha,\beta|\epsilon)=
<U_\alpha(z)V_\beta(w)>|_{z=i\epsilon;w=0}
\nonumber \\
=\sum_{{\lbrack}N|n_1...n_k\rbrack=0}^\infty\sum_{k=0}^N
{{\alpha^N\beta^k\lambda(N|k)}\over{(w-z)^{N+k}}}
=\prod_{n=1}^\infty{1\over{1-{\tilde{\alpha^n}}{\tilde{\beta}}}}
\nonumber \\
{\tilde{\alpha}}={{\alpha}\over{w-z}};
{\tilde{\beta}}={{\beta}\over{w-z}}
\nonumber
\end{eqnarray}

i.e. $G$ is the generating function for
restricted partitions with

\begin{eqnarray}
\lambda(N|k)={{(-i\epsilon)^{N+k}}\over{N!k!}}
\partial_\alpha^{N}\partial_\beta^k{G(\alpha,\beta|\epsilon)}|_{\alpha,\beta=0}
\nonumber
\end{eqnarray}

Now that we have identified the correlator generating
$\lambda(N|k)$, the next step is to identify the suitable conformal
transformation.
Note that the operator $V_\beta$ is the special
case of rank one irregular vertex operator \cite{RP},
creating a simultaneous  eigenstate of Virasoro generators $L_1$ and $L_2$
(with eigenvalues 0 and ${{\beta^2}\over2}$ respectively)
and physically can be understood as a dipole with the size $\beta$.
For this reason, it is natural to choose the transformation such that
the dipole's size shrinks to zero in the new coordinates.
So we will consider the conformal transformations of the form
{
\begin{eqnarray}
z\rightarrow{f(z)}={h(z)}e^{-{i\over{z}}}
\end{eqnarray}}
where $h(z)$ is regular and at 0, $h(0)\neq{0}$ and 
it it is smooth and analytic
in the upper half-plane (perhaps except for infinity)
 In particular,
it is instructive to consider $h(z)=1$ and $h(z)=cos{(z)}$.
Now we have to:

1. Compute infinitezimal
transformations of $U_\alpha$ and $V_\beta$.

2. Integrate them to get the finite transformations
for $U_\alpha$ and $V_\beta$ under $f(z)$.

3. Since $f(z)$ is not a fractional-linear transformation,
and its Schwarzian is singular  at 0, to match the correlators
in different coordinates,
one has to take into account the ``spontaneous symmetry breaking'' 
of the conformal symmetry (with higher Virasoro modes playing
the role of ``Goldstone modes''), by
integrating the Ward identities for $f(z)$ and regularizing the final
expression, in order
to ensure that the correlators computed in two coordinates match upon $f(z)$.

\section{ \bf Partition-Counting Correlator: the Conformal Transformations}

An important building block in 
our computation involves the finite conformal transformation 
laws for the operators of the form $T^{(n_1,n_2)}=
{1\over{n_1!n_2!}}:\partial^{n_1}\phi\partial^{n_2}\phi:$ of conformal dimensions
$n_1+n_2$ - in fact, the final answer for the number of partitions will
be expessed in terms of the finmite series in the higher-derivative Schwarzians
of $f(z)=e^{-{i\over{z}}}$.

In case of $n_1=n_2=1$ the operator $T^{(1,1)}$
(up to normalization constant
of ${1\over2}$) is just the stress-energy tensor, and both
its infinitezimal and finite transformation laws are well-known.

The infinitezimal transformation of $T^{(1,1)}$ under $z\rightarrow{z+\epsilon(z)}$ is

\begin{eqnarray}
\delta_\epsilon{T^{(1,1)}}(z)=\epsilon\partial{T^{1,1}}(z)+
2\partial\epsilon{T^{(1,1)}}(z)+{1\over{6}}\partial^3{\epsilon}(z)
\end{eqnarray}
(here and everywhere below  the infinitezimal conformal transformation parameter $\epsilon(u)$ is not
to be confused with the $i\epsilon$ for the location of $z$ which hopefully
will always be clear from the context; in our notations
the former always will appear with the argument, while the latter will not).
This infinitezimal transformation
 can be integrated to give the finite conformal transformation law 
for any $f(z)$ according to
\begin{eqnarray}
T^{(1,1)}(z)\rightarrow({{df(z)}\over{dz}})^2T^{(1,1)}(f(z))+S_{1|1}(f;z)
\end{eqnarray}
where $S_{1|1}$ is (up to the conventional normalization 
factor of ${1\over{6}}$)
the Schwarzian derivative, defined according to:
\begin{eqnarray}
S_{1|1}(f;z)= {1\over6}({{f^{\prime\prime}(z)}\over{{f^{\prime}(z)}}})^\prime
-{1\over{12}}({{f^{\prime\prime}(z)}\over{{f^{\prime}(z)}}})^2
\end{eqnarray}
The integrated transformation (3.3) can be obtained from (3.2) 
e.g. by requiring that
(3.2) is reproduced from (3.3) in the infinitezimal limit and that
the composition of two transformations $f(z)$ and $g(z)$ gives again  the
conformal transformation with $f(g(z))$.
Likewise, we can derive the transformation rules for an arbitrary
$T^{(n_1,n_2)}$.
The infinitezimal transformation is
\begin{eqnarray}
\delta_{\epsilon}T^{(n_1,n_1)}(z)=
{\lbrack}{1\over2}\oint{{du}\over{2i\pi}}
\epsilon(u):(\partial\phi)^2:(u);T^{(n_1,n_1)}(z)\rbrack
\nonumber \\
={1\over{n_2!}}:\partial^{n_1}(\epsilon\partial\phi)\partial^{n_2}\phi:(z)
+{1\over{n_1!}}:\partial^{n_1}\phi\partial^{n_2}(\epsilon\partial\phi):(z)
+{1\over{(n_1+n_2+1)!}}\partial^{n_1+n_2+1}\epsilon(z)
\end{eqnarray}
and can be integrated, by imposing the similar requirements, to give:
\begin{eqnarray}
T^{(n_1,n_2)}(z)\rightarrow{1\over{n_1!n_2!}}
:\delta_{n_1}(\phi;f(z))\delta_{n_2}(\phi;f(z)):+S_{n_1|n_2}(f;z)
\end{eqnarray}
where the operators
\begin{eqnarray}
\delta_n(\phi,f)=
{{\partial^{n}f}\over{\partial{z^n}}}\partial\phi(f(z))
+\sum_{k=1}^{n-1}\sum_{l=1}^k{{(n-1)!}\over{(n-1-k)!}}
{{\partial^{n-k}f}\over{\partial{z^{n-k}}}}
B_{k|l}(f(z);z)\partial^{l+1}\phi(f(z))
\nonumber \\
=\sum_{k=1}^n{B_{n|k}(f(z);z)}\partial^k{\phi}
\end{eqnarray}
are defined by the conformal transformation 
for $\partial^n\phi$. Here $B_{n|k}$ are the incomplete Bell polynomials
in the derivatives (expansion coefficients) of $f$
(see (3.8) and the formula below for the explicit definition) and $S_{n_1|n_2}(f;z)$
are the generalized higher-derivative Schwarzians.
To calculate $S_{n_1|n_2}$
we cast the normal ordering according to:
\begin{eqnarray}
{1\over{n_1!n_2!}}:\partial^{n_1}\phi\partial^{n_2}\phi:(z)
=lim_{\epsilon\rightarrow{0}}\lbrace
\partial^{n_1}\phi(z+{{\epsilon}\over2})\partial^{n_2}\phi
(z-{{\epsilon}\over2})+{{(-1)^{n_1}(n_1+n_2-1)!}\over{n_1!n_2!\epsilon^{n_1+n_2}}}
\rbrace
\nonumber
\end{eqnarray}
(again, the regularization parameter $\epsilon$ here is not to be confused with $i\epsilon$
for the location of $U_\alpha$ and/or infinitezimal transformation parameter $\epsilon(z)$)
Under the conformal map $z\rightarrow{f(z)}$, this expression transforms according to
\begin{eqnarray}
{1\over{n_1!n_2!}}:\partial^{n_1}\phi\partial^{n_2}\phi:(z)
\rightarrow{1\over{n_1!n_2!}}lim_{\epsilon\rightarrow{0}}\lbrace
\sum_{k_1=1}^{n_1}\sum_{k_2=1}^{n_2}{{B_{n_1|k_1}}(f(z+{{\epsilon}\over2});z)
{B_{n_k|k_2}}(f(z-{{\epsilon}\over2});z)}
\nonumber \\
\times
(\partial^{k_1}\phi(z+{{\epsilon}\over2})\partial^{k_2}\phi
(z-{{\epsilon}\over2})+{{{(-1)^{k_1}}(k_1+k_2-1)!}
\over{(f(z+{{\epsilon}\over2})-f(z-{{\epsilon}\over2}))^{k_1+k_2}}}
+{{(-1)^{n_1}(n_1+n_2-1)!}\over{n_1!n_2!\epsilon^{n_1+n_2}}})
\rbrace
\nonumber
\end{eqnarray}

Expanding in $\epsilon$, we extract (upon cancellations of the divergent 
terms) the higher-order Schwarzians to be given by:

\begin{eqnarray}
S_{n_1|n_2}(f;z)=
{1\over{n_1!n_2!}}\sum_{k_1=1}^{n_1}\sum_{k_2=1}^{n_2}
\sum_{m_1\geq{0}}\sum_{m_2\geq{0}}\sum_{p\geq{0}}
\sum_{q=1}^p
(-1)^{k_1+m_2+q}2^{-m_1-m_2}(k_1+k_2-1)!
\nonumber \\
\times
{{\partial^{m_1}B_{n_1|k_1}(f(z);z)\partial^{m_2}B_{n_2|k_2}(f(z);z)
B_{p|q}(g_1,...,g_{p-q+1})}
\over{m_1!m_2!p!(f^\prime(z))^{k_1+k_2}}}
\nonumber \\
g_s=2^{-s-1}(1+(-1)^s){{{{d^{s+1}f}\over{dz^{s+1}}}}\over{(s+1)f^\prime(z)}};
s=1,...,p-q+1
\end{eqnarray}
with the sum over the non-negative numbers $m_1,m_2$ and $p$ taken over all the combinations satisfying
$$m_1+m_2+p=k_1+k_2$$.

Here, in (3.6) and (3.7)
$B_{n|k}(g_1,...g_{n-k+1})$ are the incomplete Bell polynomials, defined 
according to:
\begin{eqnarray}
B_{n|k}(g_1,...g_{n-k+1})=n!\sum{1\over{p_1!...p_{n-k+1}!}}
{({{g_1}\over{1!}})^{p_1}}...
{({{g_{n-k+1}}\over{(n-k+1)!}})^{p_{n-k+1}}}
\end{eqnarray}
with the sum taken over all the non-negative $p_1,...p_{n-k+1}$
satisfying
\begin{eqnarray}
p_1+...+p_{n-k+1}=k
\nonumber \\
p_1+2p_2+...+(n-k+1)p_{n-k+1}=n
\nonumber
\end{eqnarray}
In particular,
the incomplete Bell polynomials $B_{n|k}(f;z)$  
 in  the derivatives 
(or the expansion coefficients) of $f(z)$,
are given by:

\begin{eqnarray}
B_{n|k}(f(z);z)=n!\sum_{n|n_1...n_k}{{\partial^{n_1}f(z)...\partial^{n_k}f(z)}
\over{n_1!...n_k!q(n_1)!...q(n_k)!}}
\nonumber
\end{eqnarray}
with the sum ${n|n_1...n_k}$ taken 
over all ordered $0<n_1\leq{n_2}...\leq{n_k}$
length $k$ partitions of $n$ and with $q(n_j)$
denoting the multiplicity of  $n_j$ element of the partition
(e.g.  for the partition $7=2+2+3$ we have $q(2)=2,q(3)=1$,
so the appropriate term would read 
$\sim{{\partial^2{f}\partial^2{f}\partial^3{f}}\over{2!2!3!\times{2!1!}}}$
).
Note that, just as the ordinary Schwarzian 
satisfies the well-known composite relation for any
combination of conformal transformations
$z\rightarrow{f(z)}\rightarrow{g(f)}$:
\begin{eqnarray}
S_{1|1}(g(z);z)=({{df}\over{dz}})^2S_{1|1}(g(f);f)+S_{1|1}(f(z);z),
\nonumber
\end{eqnarray}
the generalized Schwarzians $S_{n_1|n_2}$ also satisfy the composite
relations
\begin{eqnarray}
S_{1|1}(g(z);z)=\sum_{k_1=1}^{n_1}\sum_{k_2=1}^{n_2}
B_{n_1|k_1}(f(z);z)B_{n_2|k_2}(f(z);z)
S_{k_1|k_2}(g(f);f)+S_{n_1|n_2}(f(z);z).
\nonumber
\end{eqnarray}
(see also \cite{loriano, matone} who considered the alternative types
of higher-order Schwarzians in a rather different context).

Now let us apply the same procedure to the irregular vertex operators
$U_\alpha$ and $V_\beta$ in the partition-counting correlator (2.3).
The straightforward computation of the infinitezimal transforms
gives:

\begin{eqnarray}
\delta_\epsilon{V_\beta}
={\lbrack}\oint{{du}\over{2i\pi}}\epsilon(u)T(u);V_\beta(z)\rbrack=
:({\beta\partial{\epsilon\partial\phi}}
+{{1\over{12}}\beta^2}\partial^3\epsilon)
V_\beta:(z)
\end{eqnarray}
and
\begin{eqnarray}
\delta_\epsilon{U_\alpha}(z)
=\lbrack{\oint{{du}\over{2i\pi}}}\epsilon(u)T(u);U_\alpha(z)\rbrack
=\sum_{n=1}^\infty\lbrace
:{{\alpha^n(\partial^n(\epsilon\partial\phi))}\over{n!
(1-{{\alpha^n\partial^n\phi}\over{n!}})}}\prod_{N=1}^\infty{1\over{{1-{{\alpha^N\partial^N\phi}\over{N!}}}}}:
\nonumber \\
+
:{{{\alpha^{2n}\partial^{2n+1}\epsilon}\over{
(2n+1)!
(1-{{\alpha^n\partial^n\phi}\over{n!}})^2}}}\prod_{N=1}^\infty{1\over{{1-{{\alpha^N\partial^N\phi}\over{N!}}}}}:
\rbrace
\nonumber \\
+\sum_{0\leq{n_1}<n_2<\infty}
:{{\alpha^{n_1+n_2}\partial^{n_1+n_2+1}\epsilon}\over{
(1-{{\alpha^{n_1}\partial^{n_1}\phi}\over{{n_1}!}})
(1-{{\alpha^{n_2}\partial^{n_2}\phi}\over{{n_2}!}})}}
\prod_{N=1}^\infty{1\over{{1-{{\alpha^N\partial^N\phi}\over{N!}}}}}:
\end{eqnarray}

Integrating these infinitezimal transformations,
we obtain  the transformations of
$U_\alpha$ and $V_\beta$ for the finite conformal
transformation $f(z)=h(z)e^{-{i\over{z}}}$. For $V_\beta$, we get

\begin{eqnarray}
V_\beta(w)|_{w=i\epsilon}\rightarrow
:e^{\beta{{{\partial}f}\over{{\partial}z}}\phi(f(z))+                           
{{\beta^2}\over{2}}S_{1|1}(f(z);z)}:
|_{f(z)=h(z)e^{-{{i}\over{z}}}}
\end{eqnarray}
To determine the transformation law for $U_\alpha$, it is convenient to cast $U_\alpha$
as
\begin{eqnarray}
U_\alpha=1+\sum_{N=1}^\infty\sum_{{\lbrace}m_i\rbrace}\alpha^N\prod_{n=1}^N{({{\partial^n\phi}\over{n!}})^{m_n}}
\end{eqnarray}
with the sum over ${\lbrace}m_{i}\rbrace;i=1...N$ being taken over all the combinations
of non-negative ${\lbrace}m_{i}\rbrace$, satisfying
\begin{eqnarray}
\sum_{n=1}^N{nm_n}=N
\nonumber
\end{eqnarray}
Now introduce the {\it exchange numbers} $\nu_{ij}\geq{0}(i,j=0,...,N)$ satisfying
\begin{eqnarray}
\nu_{00}=0 \nonumber \\
\nu_{ij}=\nu_{ji} \nonumber \\
\nu_{jj}+\sum_{i=0}^N{\nu_{ij}}=m_j
\end{eqnarray}
in order to parametrize the internal normal ordering procedure for $U_\alpha$ as follows:

1.
$\nu_{ij}=\nu_{ji}(i,j\neq{0})$ defines the number of internal couplings between $(\partial^i\phi)^{m_i}$
and $\partial^j\phi^{m_j}$ factors, creating internal singularities prior to the normal ordering;

2.$\nu_{ii}$ defines the number of intrinsic same-derivative couplings between
$\partial^i\phi$'s inside each factor $(\partial^i\phi)^{m_i}$.

3. $\nu_{i0}$ counts the numbers of $\partial^i\phi$-operators left inside $(\partial^i\phi)^{m_i}$-block, that do not
participate in the contractions.

Since each coupling between $\partial^i\phi$ and $\partial^j\phi$ contributes the factor
$i!j!S_{i|j}(f;z)$ to the transformation law under $f(z)$, the overall transformation law for $U_\alpha$
is

\begin{eqnarray}
U_\alpha(z)=\prod_n:{1\over{1-{{\alpha^n\partial^n\phi}\over{n!}}}}:\rightarrow
\nonumber \\
1+\sum_{N=1}^\infty\sum_{q=1}^N\sum_{m_1,...,m_q;\lbrace\nu_{ij}\rbrace}\alpha^N
\prod_{p=1}^q{{(\delta_{p}(\phi;f))^{\nu_{p0}}m_p!(S_{p|p})^{\nu_{pp}}}\over{\nu_{p_0}!(2\nu_{pp})!!}}
\prod_{1\leq{i}<j\leq{q}}{{(S_{i|j})^{\nu_{ij}}}\over{\nu_{ij}!}}
\end{eqnarray}
where
\begin{eqnarray}
\delta_{p}(\phi;f)=\sum_{k=1}^pB_{p|k}(f(z);z)\partial^k\phi(z)
\end{eqnarray}

Finally, since the Schwarzian of the conformal transformation $f(z)$
iz nonzero,
we need to account for the spontaneous breaking of the conformal symmetry
by
 integrating the Ward identities,
in order to match the
partition-counting correlators $<U_\alpha{V_\beta}>$ in different coordinates.
For that, we first have to
 integrate the infinitezimal ``overlap''  deformation of the
correlator, emerging from the contraction of one of $\partial\phi$'s
in the stress-energy tensor with $U_\alpha$ and another with $V_\beta$.
The infinitezimal overlap deformation is given by the integral:
\begin{eqnarray}
\delta_{overlap}<U_\alpha(z)V_\beta(w)>|_{z=i\epsilon;w=0}
=\sum_{N=1}^\infty\oint{{du}\over{2i\pi}}{{\epsilon(u)}\over{(u-z)^{N+1}
(u-w)^2}}
\nonumber \\
\times:{{\alpha^N\beta}\over{(1-{{\alpha^N\partial^N\phi}\over{N!}})^2}}
\prod_{n=1;n\neq{N}}^\infty{{1}\over{1-{{\alpha^n\partial^n\phi}\over{n!}}}}:
(z):e^{\beta\partial\phi}:(w)|_{z=i\epsilon;w=0}
\nonumber \\
=
\sum_{N=1}^\infty(\partial_z^N{\lbrack}{{\epsilon(z)}\over{(z-w)^2}}{\rbrack}+
(-1)^{N+1}\partial_w{\lbrack}{{\epsilon(w)}\over{(z-w)^{N+1}}}\rbrack)
\nonumber \\
\times:{{\alpha^N\beta}\over{(1-{{\alpha^N\partial^N\phi}\over{N!}})^2}}
\prod_{n=1;n\neq{N}}^\infty{{1}\over{1-{{\alpha^n\partial^n\phi}\over{n!}}}}:
(z):e^{\beta\partial\phi}:(w)|_{z=i\epsilon;w=0}
\end{eqnarray}

This infinitezimal deformation is straightforward to integrate for the class
of the conformal transformations (2.2).
The overall integrated transformation for $U_\alpha$ under $f(z)$, with the
 overlap deformation included, is then given by
\begin{eqnarray}
U_\alpha(z)=\prod_n:{1\over{1-{{\alpha^n\partial^n\phi}\over{n!}}}}:\rightarrow
\nonumber \\
(1+\sum_{N=1}^\infty\sum_{q=1}^N\sum_{m_1,...,m_q;\lbrace\nu_{ij}\rbrace}\alpha^N
\prod_{p=1}^q{{(\delta_{p}(\phi;f))^{\nu_{p0}}m_p!(S_{p|p})^{\nu_{pp}}}\over{\nu_{p_0}!(2\nu_{pp})!!}}
\prod_{1\leq{i}<j\leq{q}}{{(S_{i|j})^{\nu_{ij}}}\over{\nu_{ij}!}})
\nonumber \\
\times
\prod_{n=1}^\infty
{1\over{1-{{\alpha^n\beta}\over{n!}}{{D_n(f(z);z)}\over{(1-{{\alpha^n\delta_n(\phi,f)}})}}}}
\end{eqnarray}
where
\begin{eqnarray}
D_n(f(z);z)
\equiv{D(n)}
=-i\sum_{k=1}^n(-1)^{k+1}k!{{B_{n|k}(f(z);z)}\over{(f(z)-f(0))^{k}}}
\end{eqnarray}
(to abbreviate notations, below we will also use the symbol $D(n)$ for
$D_n(f(z);z)$) 
For the conformal transformations of the form
$f(z)=h(z)e^{-{i\over{z}}}$ the overall transformation law for $V_\beta(z)$ remains the same
up to  terms that vanish identically at $z=0$:
\begin{eqnarray}
V_\beta(w)|_{w=i0}\rightarrow
e^{\beta{{{d}f}\over{{d}z}}
\partial\phi(f(z))+                  
{{\beta^2}\over{2}}S_{1|1}(f(z);z)}
|_{f(z)=h(z)e^{-{{i}\over{z}}}}
\end{eqnarray}
Now that we are prepared to calculate
the partition-counting correlator in the new coordinates, here
comes the crucial part.
The dipole's size in the new coordinates is
\begin{eqnarray}
\beta{{df}\over{dz}}|_{z{\rightarrow}0}\rightarrow{0}
\end{eqnarray}
and shrinks to zero with our choice of $f(z)$.
This drastically simplifies the calculation.
While the operator $U_\alpha$ looks extremely cumbersome
in the new coordinates (particularly, because of the complexities
involving $\delta_n(\phi;f)$-operators), any contractions of
derivatives of $\phi$ in $U_\alpha$ with $V_\beta$ bring down
the factors proportional to $f^\prime(z)$ and therefore 
vanish for the conformal transformations of the form (2.2).
As a result, only the zero modes of $U_\alpha$ and $V_\beta$ contribute 
to the correlator in the new coordinates.
Technically, this implies $\nu_{0j}=0$ for all $j$.
The correlator is then easily computed to give the generating function
for the restricted partitions in terms of higher-derivative Schwarzians and incomplete Bell polynomials:
\begin{eqnarray}
G(\alpha,\beta|\epsilon){\equiv}<U_\alpha(z)V_\beta(w)>|_{z=i\epsilon,w=0}
=e^{{1\over2}\beta^2S_{1|1}(f(w);w)}|_{w=0}
\nonumber \\
\times
(1+\sum_{N=1}^\infty\sum_{q=1}^N\sum_{m_1,...,m_q;\lbrace\nu_{ij}\rbrace}\alpha^N
\prod_{p=1}^q{{m_p!(S_{p|p})^{\nu_{pp}}}\over{(2\nu_{pp})!!}}
\prod_{1\leq{i}<j\leq{q}}{{(S_{i|j})^{\nu_{ij}}}\over{\nu_{ij}!}})
\times
\prod_{n=1}^\infty
{1\over{1-{{\alpha^n\beta}\over{n!}}{{D(n)}}}}|_{z=i\epsilon}
\end{eqnarray}
The overall constant ($\epsilon$-independent)
factor of $e^{{1\over2}\beta^2S_{1|1}(f;w)}|_{w=0}$ is related
to the Casimir energy associated with the conformal transformation $f(w)$.
It is  irrelevant and disappears 
when the correlator is normalized with the inverse
of the partition function of the system. The generation function for
the partition numbers is then simply obtained 
by replacing this factor with 1.

Now the final step is to take the derivatives of $G$  in $\alpha,\beta$
and to $\epsilon$-order  the result, retaining the finite terms as
$\epsilon$ is set to 0.
Straightforward calculation gives:
\begin{eqnarray}
\lambda(N|Q)=
(-i\epsilon)^{N+Q}:\sum_{N_1=0}^{N}\sum_{N_2=0}^{N-N_1}
\sum_{\lbrace{m_j\geq{0}}\rbrace}\sum_{\lbrace{n_j;p_j{\geq}0}\rbrace}
\sum_{\lbrace{\nu_{ij}\geq{0}};\rbrace}
\nonumber \\
{{\prod_{k=1}^{N_1}\prod_{q,r;1\leq{q}<r\leq{N_1}}:{{m_k!(S_{k|k})^{\nu_{kk}}(S_{q|r})^{\nu_{qr}}}\over{
(2\nu_{kk})!\nu_{qr}!}}}}
(D(n_1))^{p_1}...(D(n_Q))^{p_Q}:_{N+Q}
\end{eqnarray}
with the summations/products taken over the non-negative integer values
of 
$$N_1,N_2,N_3;m_1,...m_{N_1};n_1,...n_Q;p_1,...p_Q,\lbrace{\nu_{ij}}\rbrace;1\leq{i,j}\leq{N_1}$$
 satisfying:
\begin{eqnarray}
N_1+N_2=N
\nonumber \\
\sum_{j=1}^{N_1}jm_j=N_1
\nonumber \\
\sum_{j=1}^sn_jp_j=N_2
\nonumber \\
\nu_{ii}+\sum_{j=1}^{N_1}\alpha_{ij}=m_i;1\leq{i}\leq{m_i}
\nonumber \\
\nu_{ij}=\nu_{ji}
\end{eqnarray}
The $\epsilon$-ordering symbol $:...:_{N+Q}$ in each monomial term of the sum $\sim:S...SD...D:_{N+Q}$
by definition only retains the terms of the order of $\epsilon^{-N-Q}$  
upon the evaluation of each product $S...SD...D$, in order to ensure that the overall contribution is finite,
upon multiplication by $(-i\epsilon)^{N+Q}$
(we refer to this procedure as $\epsilon$-ordering to distinguish it from the usual normal ordering
defined for operators in CFT).
Let us stress that, since each $S$ or $D$ has the finite and definite singularity order in $\epsilon$,
the overall result for $\lambda(N|Q)$ is the exact analytic expression,
given by the finite series, uniquely determined  by the structures of $S$ and $D$ for each $f$
(with $f$ satisfying the constraints described above).  
This concludes our derivation of counting the restricted partitions,
expressed in terms of finite series in the incomplete Bell polynomials
and the generalized higher-derivative Schwarzians of the
defining conformal transformation $f(z)=h(z)e^{-{{i}\over{z}}}$.

\section{\bf Conclusion. Tests and comments}

Having presented the exact analytic expressionb for the number of the 
partitions, in this section we will
provide some checks and examples of how the expression (3.22), 
constituting the main result of this paper,
 works in 
practice.
First of all, it is quite straightforward
 to demonstrate that the expression (3.22)
leads to the correct answer for any partition number in
the case of the conformal transformation $f(z)$ with
the simplest choice $h(z)=1$. 
Let us start from the most elementary case of 
the maximal length partition where obviously $\lambda(N|N)=1$.
for any $N$. Indeed, according to (3.22),
one has
\begin{eqnarray}
\lambda(N|N)=(-i\epsilon)^{2N}:(D(1))^N:_{2N}
=(-i\epsilon)^{2N}(-i)^N(-{{i\over{\epsilon^2}}})^N=1
\end{eqnarray}
Similarly, it is easy to verify
the case $1=\lambda(N-1|N)$
Indeed,
\begin{eqnarray}
\lambda(N|N-1)=
(-i\epsilon)^{2N-1}:(D(1))^{N-2}D(2):_{N-1}
\nonumber \\
=:(-i\epsilon)^{2N-1}(-{{i\over{\epsilon}}})^{2N-4}({i\over{\epsilon^3}}
-{1\over{2\epsilon^4}})+O(\epsilon):_{2N-1}=1
\end{eqnarray}
Note ( although irrelevant to our result) the appearance
of the singular term $\sim{1\over{\epsilon}}$ at this level
(which was absent in the case of $\lambda(N|N)$). This term 
disappears upon the $\epsilon$-ordering procedure and is of no significance for our purposes,
 but the very reason
for its emergence
is also related to the Schwarzian singularities at $\epsilon\rightarrow{0}$.
It is easy to check that the results (4.1), (4.2) 
actually hold for any smooth regular $h(z)$ satisfying the conditions
defined above, not just for $h=1$. However, the case of $h=1$ is the easiest one to verify the correctness
of (3.22) for any partition, as this can be done by simple analysis of the $\epsilon$-dependence.
Indeed, in general, each term for the partition number $\Lambda(N|Q)$ in (3.22) typically consists
of $Q$ $D(n)$-factors and R $S$-factors (Schwarzians of all kinds and orders), where $R$ can in principle vary as
$0\leq{R}\leq{{\lbrack{{N-Q}\over2}}}\rbrack$
However, in case of $h(z)=1$ only the terms with $R=0$ contribute. Indeed, the terms in each of the 
Schwarzians $S_{n_1|n_2}$, least singular in $\epsilon$ , are of the order of ${1\over{\epsilon^{n_1+n_2+2}}}$
(e.g. $S_{1|1}(i\epsilon)={1\over{12\epsilon^4}}$.  On the other hand, the $D(n)$-factors, consisting of
combinations of incomplete Bell polynomials $B_{n|k}(f(z);z)$  with various $k$'s,  have the lowest
singularity order $\sim{1\over{\epsilon^{n+1}}}$ for $h=1$. Thus it is clear that each contribution with nonzero $R$ ,
upon multiplication by $(i\epsilon)^{N+Q}$ has the  singularity order of at least ${1\over{\epsilon^R}}$ and will
disappear upon the normal ordering $:...:_{N+Q}$. Furthermore, 
the only source of  the lowest singularity terms in $D(n)$ of the order of ${1\over{\epsilon^{n+1}}}$
is $B_{n|1}$, as all other $B_{n|k}$ with $k>1$ are more singular, as it is easy to check. 
These terms stem from the derivatives $\partial^n(e^{-{i\over{z}}})|_{z=i\epsilon}$  and are easily 
computed to be given by (skipping terms with the higher order singularities, not contributing
to the $\epsilon$-ordering)
$$-{i\over{n!}}(e^{{i\over{z}}})\partial^n(e^{-{i\over{z}}})|_{z=i\epsilon}=({i\over{\epsilon}})^{n+1}
+h.s.$$
Thus each of the terms with $R=0$:
$$\sim(-i\epsilon)^{N+Q}:D(n_1)...D(n_Q):_{N+Q}{(N=n_1+...+n_Q)}=(-i\epsilon)^{N+Q}({i\over{\epsilon^{N+Q}}})=1$$
contributes 1 
to the sum. But the number of such terms obviously equals the number of partitions
$\lambda(N|Q)$, hence this constitutes the proof that the formula (3.22) 
works correctly with the conformal transformation
$f(z)=e^{-{i\over{z}}};h(z)=1$. Although the case of $h(z)=1$ is somewhat simplistic
(e.g. with no Schwarzians entering the game), this by itself is already a non-trivial check of how
the conformal invariance works in (3.22).
Of course, with $h(z)\neq{1}$ things change significantly and the Schwarzians of all orders contribute 
nontrivially
to the expression (3.22) for the partitions.
For example, consider $h(z)=cos(z)$ and $\lambda(4|2)=2$.
In this case, the Schwarzian $S_{1|1}$ is given by
\begin{eqnarray}
6S_{1|1}|_{z=i\epsilon}={1\over{2z^4}}+{{2i}\over{z}}={1\over{2\epsilon^4}}+{2\over{\epsilon}}+O(\epsilon)
\end{eqnarray}
and does of course contribute to the normal ordering in general (the same is true for other 
$S_{n_1|n_2}$'s).
According to (3.22) we have
\begin{eqnarray}
\lambda(4|2)=(-i\epsilon)^6{\lbrack}:S_{1|1}(D(1))^2:_{6}+:(D(2))^2:_{6}+:D(1)D(3):_{6}\rbrack
\end{eqnarray}
Straightforward calculation gives:
\begin{eqnarray}
(-i\epsilon)^6:S_{1|1}(D(1))^2:_6=0 \nonumber \\
(-i\epsilon)^6:(D(2))^2:_6=1 \nonumber \\
(-i\epsilon)^6:D(3)D(1):_6=1\nonumber \\
\lambda(4|2)=2
\end{eqnarray}
For this partition, the Schwarzian related term still does not contribute, although
its vanishing is not automatic but is related to the particular $\epsilon$-structure of the Schwarzian
$S_{1|1}(f(z);z)$ for $h(z)=cos(z)$.
For $\lambda(5|2)$ one calculates:
\begin{eqnarray}
\lambda(5|2)=(-i\epsilon)^7{\lbrack}
:S_{1|2}(D(1))^2:_{7}+:S_{1|1}D(1)D(2):_{7}+:D(1)D(4):_{7}+:D(2)D(3):_{7}\rbrack\nonumber \\
(-i\epsilon)^7:S_{1|1}D(1)D(2):_7={7\over{12}}\nonumber \\
(-i\epsilon)^7:S_{1|2}(D(1))^2:_7=-{1\over4} \nonumber \\
(-i\epsilon)^7:D(1)D(4):_7={1\over{4}} \nonumber \\
(-i\epsilon)^7:D(2)D(3):_7={{17}\over{12}} \nonumber \\
\lambda(5|2)=2
\end{eqnarray}
so for this partition both $S$-type and $D$-type terms  contribute nontrivially to
$\lambda$. One can perform  some similar tests to show that (3.22) works correctly.
In general, however, the complexity of the manifest expressions for $\lambda(N|Q)$ grows dramatically
with $N$ and especially with the difference $N-Q$, as not only the structure of higher order
Schwarzians becomes increasingly cumbersome, but also the $\epsilon$-ordering procedure of the 
terms gets quite tedious. For this reason, the formula (3.22), although exact, is in practice
 less convenient
for numerical computations of the partitions, compared to using the standard generating functions.
Nevertheless, it casts the partition numbers in terms of exact finite analytic expressions
in terms of the conformal transformation (2.2),
which demonstrates the power of conformal symmetry and constitutes the main result of this work.

\begin{center}
{\bf Acknowledgements}
\end{center}

The author acknowledges the support of this work by  the National Natural 
Science Foundation of China under grant 11575119.
I also would like to express my gratitude to Hermann Nicolai and Rakibur Rahman
for their kind hospitality at Max Planck Institute for Gravitational Physics 
(Albert Einstein Institute) in Potsdam, where the concluding part of this work has been
done.


\end{document}